\documentclass[11pt]{article}

\usepackage{amssymb}
\usepackage{amscd}
\usepackage{amsmath}
\usepackage{amsfonts}
\usepackage{theorem}
\usepackage{mathrsfs}

{\theorembodyfont{\slshape}
\newtheorem{theorem}{Theorem}

\newtheorem{lemma}{Lemma}
\newtheorem{corollary}{Corollary}

}
{\theorembodyfont{\rmfamily}
\newtheorem{definition}{Definition}

\newtheorem{remark}{Remark}

}
\def\proof{{\noindent\sc Proof. \quad}}
\newcommand{\proofof}[1]{{\noindent\sc Proof of #1. \quad}}
\def\eproof{{\mbox{}\hfill\qed}\medskip}
\newcommand\qed{{\unskip\nobreak\hfil\penalty50\hskip2em\vadjust{}
\nobreak\hfil$\Box$\parfillskip=0pt\finalhyphendemerits=0\par}}

\def\R{{\mathbb{R}}}
\def\E{{\mathbb{E}}}
\def\Prob{{\rm{Prob}}}

\def\K{{\cal{K}}}
\def\U{{\cal{U}}}
\def\B{{\cal{B}}}
\def\S{{\cal{S}}}
\def\D{{\cal{D}}}
\def\Oh{{\cal{O}}}
\def\msD{\mathscr D}
\def\st{{\rm{s.t.}}}
\def\transp{^{\rm T}}
\def\mtransp{^{-{\rm T}}}
\def\sign{\mathsf{sign}}
\def\sfd{\mathsf{d}}
\def\Id{\mathsf{Id}}
\def\Sing{\mathsf{Sing}}
\def\uno{\mbox{1\hspace*{-2.5pt}l}}
\def\u{{\mathsf{u}}}
\def\fG{{\mathfrak{G}}}

\begin{document}

\bibliographystyle{plain}

\makeatletter


\def\JACM{Journal of the ACM}
\def\CACM{Communications of the ACM}
\def\ICALP{International Colloquium on Automata, Languages
            and Programming}
\def\STOC{annual ACM Symp. on the Theory
          of Computing}
\def\FOCS{annual IEEE Symp. on Foundations of Computer Science}
\def\SIAM{SIAM Journal on Computing}
\def\SIOPT{SIAM Journal on Optimization}
\def\MOR{Math. Oper. Res.}
\def\BSMF{Bulletin de la Soci\'et\'e Ma\-th\'e\-ma\-tique de France}
\def\CRAS{C. R. Acad. Sci. Paris}
\def\IPL{Information Processing Letters}
\def\TCS{Theoretical Computer Science}
\def\BAMS{Bulletin of the Amer. Math. Soc.}
\def\TAMS{Transactions of the Amer. Math. Soc.}
\def\PAMS{Proceedings of the Amer. Math. Soc.}
\def\JAMS{Journal of the Amer. Math. Soc.}
\def\LNM{Lect. Notes in Math.}
\def\LNCS{Lect. Notes in Comp. Sci.}
\def\JSL{Journal for Symbolic Logic}
\def\JSC{Journal of Symbolic Computation}
\def\JCSS{J. Comput. System Sci.}
\def\JoC{J. of Complexity}
\def\MP{Math. Program.}
\sloppy

\begin{title}
{{\bf  \mbox{Solving Linear Programs with Finite Precision:} \\
III. Sharp Expectation bounds}
\thanks{This work has been substantially funded by
a grant from the Research Grants Council of the
Hong Kong SAR (project number CityU 1085/02P).
}}
\end{title}
\author{Dennis Cheung\\
Division of Continuing Professional Education\\
The Hong Kong Institute of Education\\ 
HONG KONG\\
e-mail: {\tt chinwing@ied.edu.hk}
\and
Felipe Cucker\\
Department of Mathematics\\
City University of Hong Kong\\
HONG KONG\\
e-mail: {\tt macucker@cityu.edu.hk}
}

\date{}

\makeatletter
\maketitle
\makeatother

\thispagestyle{empty}
 
\begin{quote}
{\small {\bf Abstract.\quad}
We give an $\Oh(\log n)$ bound for the expectation of the 
logarithm of the condition number $\K(A,b,c)$ introduced in 
``Solving linear programs with finite 
precision: I. Condition numbers and random programs.''
\mbox{\it Math. Programm.}, 99:175--196, 2004. This bound 
improves the previously existing bound, which was of 
$\Oh(n)$.  
}
\end{quote}

\section{Introduction}

Consider the following linear programming problem 
(in standard form), 
\begin{align}\label{K3:LP}
  \min c\transp x&\notag\\
  \st \ Ax &=b \tag*{(P)}\\
     x& \geq 0.\notag
\end{align}
Here $A \in \R^{m\times n}, b \in\R^m, c \in\R^n$, 
and $n\geq m\geq 1$.

Assuming this problem is feasible (i.e., the set given by 
$Ax=b$, $x\geq 0$, is not empty) and bounded (i.e., 
the function $x\mapsto c\transp x$ is bounded below 
on the feasible set), algorithms solving \ref{K3:LP}  
may return an optimizer $x^*\in\R^n$ and/or 
the optimal value $c\transp x^*$. Whereas these two 
computations are essentially equivalent in the presence 
of infinite precision, obtaining an optimizer appears to 
be more difficult if only finite precision is available. 
Accuracy analyses of interior-point algorithms for 
these problems have been done in~\cite{Vera98} 
---for the computation of the optimal value--- 
and in~\cite{ChC03} ---for the computation of 
an optimizer. In both cases, accuracy  bounds 
(as well as complexity bounds) are given in terms 
of the dimensions $m$ and $n$, as well as of the 
logarithm of a condition number. The bounds in 
both analyses are similar. What turns out to be 
different is their relevant condition numbers. 

In~\cite{Vera98} this is Renegar's condition number 
$C(A,b,c)$ which, 
roughly speaking, is the relativized inverse of the 
size of the smallest perturbation needed to 
make~\ref{K3:LP} either infeasible or unbounded. 
In~\cite{ChC03} it is the condition number 
$\K(A,b,c)$ which, following the same idea, 
is the relativized inverse of the 
size of the smallest perturbation needed to 
change the optimal basis of~\ref{K3:LP} 
(a detailed definition is in Section~\ref{K3:sec:prelim}
below). 

A characteristic of these (and practicality all other) condition numbers 
is that they cannot be easily computed from the data at hand. Their 
computation appears to be at least as difficult as that of the solution for 
the problem whose condition they are 
measuring (see~\cite{Renegar94} for a discussion on this)  
and requires at 
least the same amount of precision (see~\cite{ChC05}). 
A way out of this dilemma going back to the very beginning 
of condition numbers is to randomize the data and to 
estimate the expectation of its condition. Indeed, 
the first papers on condition are published independently 
by Turing~\cite{Turing48} and by Goldstine and 
von Neumann~\cite{vNGo47}, both for the condition of 
linear equation solving and in a sequel~\cite{vNGo51} to the 
latter the matix $A$ of the input linear system was considered 
to be random and some probabilistic estimates on its 
condition number were derived. This approach was 
subsequently championed 
by Demmel~\cite{Demmel88} and Smale~\cite{Smale97}. 

A number of probabilistic estimates for Renegar's condition number 
(or for a close relative introduced in~\cite{ChC00}) have been 
obtained in the last 
decade~\cite{CW01,BuCuLo:AoP,DuSpTe:09}.  
The overall picture 
is that the contribution of the log of this condition number to complexity 
and accuracy bounds is, on the average, $\Oh(\log n)$. 
In contrast with this satisfactory state of affairs, 
little is known for the condition number $\K$ on random 
triples $(A,b,c)$. In~\cite{ChC02} it was shown that 
for these triples, conditioned to~\ref{K3:LP} being feasible and 
bounded, $\log \K(A,b,c)$ is $\Oh(n)$ on the average but this 
estimate appears to be poor. In the present paper we improve 
this result and show a $\Oh(\log n)$ bound (see 
Theorem~\ref{K3:th:mainK} below for a precise statement).

\section{Statement of the Main Result}
\label{K3:sec:prelim}

In this section we fix notations, recall the definition 
of $\K(A,b,c)$, and state our main result. 

For any subset $B$ of $\{1, 2, . . . , n\}$, denote by $A_B$ the
submatrix of $A$ obtained by removing from $A$ all the 
columns with index not in $B$. If $x\in \R^n$, 
$x_B$ is defined analogously.
A set $B\subset\{1, 2, . . . , n\}$ such that $|B| = m$ and $A_B$ is
invertible is said to be a {\em basis} for $A$.

Let $B$ be a basis. Then we may uniquely solve $A_Bx'  = b$.
Consider the point $x^*\in \R^n$ defined by $x^*_j = 0$ for
$j\not\in B$ and $x^*_B = x'$. Clearly, $Ax^* = b$. We say that
$x^*$ is a {\em primal basic solution}. If, in addition, $x^*\geq 0$,
which is equivalent to $x^*_B\geq 0$, then we say 
$x^*$ is a {\em primal basic feasible solution}.

The dual of \ref{K3:LP}, which in the sequel we denote by \ref{K3:LD}, 
is the following problem,
\begin{align}\label{K3:LD}
\max\ b\transp y&\tag*{(D)}\\
\st\ A\transp y&\leq c.\notag
\end{align}

For any basis $B$, we may now uniquely solve 
$A\transp _By^* = c_B$. The point $y^*$ thus 
obtained is said to be a {\em dual basic solution}. If,
in addition, $A\transp y \leq c$, $y^*$ is said to be a 
{\em dual basic feasible solution}.

Let $B$ be a basis. We say that $B$ is an 
{\em optimal basis} (for the pair (P--D)) 
if both the primal and dual basic solutions are
feasible. In this case the points $x^*$ and $y^*$ 
above are the {\em optimizers} of~\ref{K3:LP} 
and~\ref{K3:LD}, respectively.   

We denote by $d$ the input data $(A, b, c)$. We say that $d$ 
is {\em feasible} when there exist $x\in\R^n$, $x\geq 0$, 
and $y\in\R^m$ such that $Ax=b$ and $A\transp y\leq c$. 
Let
$$
  \U=\{d=(A,b,c)\mid\mbox{$d$ has a unique optimal basis}\}.
$$
By definition, triples in $\U$ are feasible. 

To define conditioning, we need a norm in the space of data triples.
To do so, we associate to each triple 
$d = (A, b, c)\in\R^{mn+m+n}$ the matrix
$$
M_d = \left(\begin{array}{cc}
  c\transp  &0\\
   A& b\\
  \end{array}\right)
$$
and we define $\|d\|$ to be the operator norm
$\|M_d\|_{rs}$ of $M_d$ considered as a linear map
from $\R^{n+1}$ to $\R^{m+1}$. Note that this 
requires norms $\|\ \|_r$ and $\|\ \|_s$ in 
$\R^{n+1}$ and $\R^{m+1}$, respectively.

Let $\Sigma_U$ be the boundary of $\U$ in 
$\R^{mn+m+n}$. For any data input $d\in\U$,
we define the {\em distance to ill-posedness} and the 
{\em condition number} for $d$, respectively, as follows,
$$
  \varrho(d) = \min\{ \|\delta d\| :\, d+\delta d\in
  \Sigma_U\}\qquad\mbox{ and }
  \qquad \K(d) = \frac{\|d\|}{\varrho(d)}.
$$

We next state our main result, after making precise the 
underlying probability model. 

\begin{definition}
We say that $d=(A,b,c)$ is {\em Gaussian}, and 
we write $d\sim N(0,\Id)$, when 
all entries of $A, b$ and $c$ are i.i.d. with standard
normal distribution.
\end{definition}

\begin{theorem}\label{K3:th:mainK}
For the $\|\ \|_{12}$ norm we have 
$$
  \underset {d\sim N(0,\Id)}{\E}
  \left(\ln\K(d)\mid d\in\U\right) \leq
  \frac{5}{4}\ln(m+1)+\frac{3}{2}\ln(n+1)+\ln(12).
$$
\end{theorem}

\begin{remark}
The use of the $\|\ \|_{12}$ norm in Theorem~\ref{K3:th:mainK} is 
convenient but inessential. Well known norm equivalences 
yield $\Oh(\log n)$ bounds for any of the usually 
considered matrix norms.
\end{remark}

\section{Proof of the Main Result}

\subsection{A useful characterization}

Write $\D=\R^{mn+m+n}$ for the space of data inputs, and 
$$
  \B = \{B\subset\{1,2,\ldots,n\}|\, |B| = m\}
$$
for the family of possible bases.

For any $B\in\B$ and any triple $d\in\D$, 
let $\S_1$ be the set of all $m$ by $m$
submatrices of $[A_B, b]$, $\S_2$ the set of all $m + 1$ by
$m + 1$ submatrices of $\left(A\transp , c \right)\transp $ 
containing $A_B$, and $\S_B(d) = \S_1 \bigcup \S_2$. 
Note that $|\S_1| = m + 1$ and
$|\S_2| = n-m$, so $\S_B(d)$ has $n + 1$ elements.

Let $\Sing$ be the set of singular matrices. For any square matrix
$S$, we define the distance to singularity as follows.
$$
  \rho_\Sing(S) := \min\{ \|\delta S\| :\, (S + \delta S) \in\Sing\}.
$$
For any $B\in\B$ consider the function
\begin{eqnarray*}
   h_B:\D &\to& [0,+\infty)\\
           d&\mapsto & \min_{S\in\S_B(d)} \rho_\Sing(S).
\end{eqnarray*}

The following characterization of $\varrho(d)$ 
is Theorem~2 in~\cite{ChC02}.

\begin{theorem}\label{K3:the1}
For any $d\in\U$,
$$
   \varrho(d) = h_B(d)
$$
where $B$ is the optimal basis of $d$.
\eproof
\end{theorem}

\subsection{The group action}

We consider the group (with respect to componentwise 
multiplication) $\fG_n=\{-1,1\}^n$. This group acts 
on $\D$ as follows. For $\u\in\fG_n$ let 
$D_{\u}$ be the diagonal matrix having $\u_j$ 
as its $j$th diagonal entry, and 
\begin{eqnarray*}
\u(A)&:=&AD_{\u} = (\u_1a_1,\u_2a_2,\ldots,
\u_{n}a_{n}),\\
  \u(c)&:=&D_{\u}c =
 (\u_1c_1,\u_2c_2,\ldots, \u_{n}c_{n}),
\end{eqnarray*} 
where $a_i$ denotes the $i$th column of $A$. 
We define $\u(d) := (\u(A), b, \u(c))$. The group $\fG_n$ also acts 
on $\R^n$ by $\u(x):=(\u_1x_1,\ldots,\u_nx_n)$. It is immediate 
to verify that for all $A\in\R^{m\times n}$, all $x\in\R^n$, and 
all $\u\in\fG_n$ we have $\u(A)\u(x)=Ax$.

\begin{lemma}\label{K3:l1}
The functions $h_B$ are $\fG_n$-invariant. That is, 
for any $d\in\D$, $B\in\B$ and $\u\in\fG_n$,
$$
   h_B(d)=h_B(\u(d)).
$$
\end{lemma}

\proof
Let $S^*$ be any matrix in $\S_B(d)$ such that
\begin{equation}\label{K3:eq0}
   \rho_\Sing(S^*) = \min_{S\in\S_B(d)}\rho_\Sing(S).
\end{equation}
Let $k$ be the number of rows (or columns) of $S^*$ and 
$E$ be any matrix in $\R^{k\times k}$ such that 
$S^*+ E\in \Sing$ and 
\begin{equation}
\|E\|=\rho_\Sing(S^*)\label{K3:eq1}.
\end{equation}
Then, there exists $z\in\R^k$ such that
\begin{equation}\label{K3:eq11}
   (S^* + E)z =0.
\end{equation}

Suppose $S^*$ consists of the $j_1, j_2, \ldots, j_k$ 
columns of $M_d$ and let 
$\bar{\u}=(\u_{j_1}, \u_{j_2},\ldots,\u_{j_k})\in\fG_k$. 
Then, by the definition of $\S_B(d)$ and $\S_B(\u(d))$, we 
have $\bar{\u}(S^*)\in\S_B(\u(d))$.
Furthermore, 
$$
 (\bar{\u}(S^*) +\bar{\u}(E))\bar{\u}(z) 
 = \bar{\u}(S^*+E)\bar{\u}(z) 
 = (S^*+E)(z) =0,
$$
the last by Equation~\eqref{K3:eq11}.
That is, $(\bar{\u}(S^*)+\bar{\u}(E))$ is also singular. By the
definition of $\rho_\Sing$,
\begin{equation}
\rho_\Sing(\bar{\u}(S^*))\leq \|\bar{\u}(E)\|.\label{K3:eq2}
\end{equation}
Since operator norms are invariant under 
multiplication of arbitrary matrix columns by $-1$ 
we have $\|E\|= \|\bar{\u}(E)\|$. 
Combining this equality with Equations (\ref{K3:eq0}), (\ref{K3:eq1}), 
and (\ref{K3:eq2}) we obtain
\begin{equation*}
  \rho_\Sing(\bar{\u}(S^*))\leq \min_{S\in\S_B(d)}\rho_\Sing(S).
\end{equation*}
Since $\bar{\u}(S^*)\in\S_B(\u(d))$ we obtain 
\begin{equation*}\label{K3:eq4}
  \min_{S\in\S_B(\u(d))}\rho_\Sing(S)\leq
  \min_{S\in\S_B(d)}\rho_\Sing(S).
\end{equation*}
The reversed inequality follows by exchanging 
the roles of $S(\u)$ and $S$. 
\eproof

For any $B\in\B$, let
$$
   \U_B=\{d\in\D\mid B \mbox{ is the only optimal basis for $d$}\}.
$$
The set $\U$ of well-posed feasible triples is thus partitioned 
by the sets $\{\U_B\mid B\in\B\}$.

\begin{lemma}\label{K3:l2}
Let $d\in\D$ and $B\in\B$. 
If $h_B(d)>0$, then there exists a unique $\u\in\fG_n$
such that $\u(d)\in\U_B$.
\end{lemma}

\proof
First observe that, since $\min_{S\in\S_B(d)}\rho_\Sing(S)>0$, 
we have $A_B$ invertible and therefore $B$ is a basis for $A$. 
Let $y^*$ and $x^*$ be the dual and primal basic solutions of
$d$ for the basis $B$, i.e.
\begin{equation}\label{K3:eq6}
   y^*=A_B\mtransp c_B,\quad x^*_B 
  =A_B^{-1}b,\quad x^*_j = 0,\, \forall j\not\in B.
\end{equation}
Similarly, let $y^{\u}$ and $x^{\u}$ be the dual and
primal basic solutions of $\u(d)$ for the same basis. Then, 
using that $\u(A)=AD_{\u}$ and $\u(c)=D_{\u}c$, 
\begin{equation}\label{K3:eq7y}
  y^{\u}=\u(A)_B\mtransp\,\u(c)_B 
  =A_B\mtransp (D_\u)_B\mtransp (D_\u)_Bc_B 
  =A_B\mtransp c_B =y^*
\end{equation}
the third equality by the definition of $(D_\u)_B$. 
Similarly, 
\begin{equation}\label{K3:eq7x}
  x^{\u}_B
  =\u(A)_B^{-1} \,b 
  =(D_\u)_B^{-1}A_B^{-1} \,b 
  =(D_\u)_B A_B^{-1}\,b
  =(D_\u)_B x^*_B
\end{equation}
and $x^{\u}_j = 0$ for all $j\not\in B$.
Therefore,
\begin{eqnarray}\label{K3:eq9}
   B \mbox { is optimal for $\u(d)$}
&\Leftrightarrow&\mbox{$x^{\u}$ and $y^{\u}$ 
     are both feasible}\notag\\
&\Leftrightarrow&\left\{\begin{array}{l} x^{\u}_B\geq 0\\[3pt]
\u(A)_j\transp y^\u\leq \u(c)_j,\ \mbox{for }j\not\in B
\end{array}\right.\notag\\
&\Leftrightarrow&\left\{\begin{array}{l} (D_\u)_Bx^*_B\geq 0\\[3pt]
\u_j(a_j)\transp y\leq \u(c)_j,\ \mbox{for }j\not\in B
\end{array}\right.\qquad\mbox{(by (\ref{K3:eq7y}) and (\ref{K3:eq7x}))}\notag\\
&\Leftrightarrow&
\left\{\begin{array}{l} \u_jx^*_j\geq 0,\ \mbox{for }j\in B\\[3pt]
\u_j(c_j-a_j\transp y)\geq 0,\ \mbox{for }j\not\in B.
\end{array}\right.
\end{eqnarray}

Since by hypothesis $\min_{S\in\S_B(d)}\rho_\Sing(S)>0 $,
\begin{equation}\label{K3:eq10}
x^*_j\neq 0,\,\forall j\in B\qquad\mbox{ and }\qquad 
a_j\transp y\neq
c_j,\,\forall j\not\in B.
\end{equation}

Combining Equations (\ref{K3:eq9}) and (\ref{K3:eq10}), 
the statement follows for $\u\in\fG_n$ given by 
$\u_j=\sign(x^*_j)$ if $j\in B$ and 
$\u_j=\sign(c_j-a_j\transp y)$ otherwise. Clearly, this 
$\u$ is unique.
\eproof

For $B\in \B$ let 
$$
  \Sigma_B:=\Big\{d\in\D\mid h_B(d)=0 \Big\}
$$
and $\D_B:=\D\setminus \Sigma_B$. Lemma~\ref{K3:l1} 
implies that, for all $B\in\B$, $\Sigma_B$ and $\D_B$ are 
$\fG_n$-invariant. Lemma~\ref{K3:l2} immediately 
implies the following 
corollary.

\begin{corollary}\label{cor:K3_partition}
For all $B\in\B$ the sets 
$$
    \D_{\u}:=\{d\in \D_B\mid \u(d)\in \U_B\},\qquad 
    \mbox{for $\u\in\fG_n$}
$$ 
are a partition of $\D_B$. \eproof
\end{corollary}


\subsection{Probabilities}

\begin{definition}\label{K3:def:invariance}
We say that a distribution $\msD$ on the set of triples $d=(A,b,c)$ 
is {\em $\fG_n$-invariant} when 
\begin{description}
\item[(i)]
   if $d\sim \msD$ then $\u(d)\sim\msD$ for all $\u\in\fG_n$.
\item[(ii)]
  for all $B\in\B$, 
 $\displaystyle    \underset {d\sim \msD}{\Prob}
   \big\{h_B(d)=0\big\}=0.$
\end{description}
\end{definition}

Note that Gaussianity is  a special case of 
$\fG_n$-invariance. Consequently, all results true for a 
$\fG_n$-invariant distribution also hold for Gaussian data.
\medskip

\noindent
{\bf Note:\ }
For a time to come we fix a $\fG_n$-invariant 
distribution $\msD$ with density function~$f$. 
\smallskip

\begin{lemma}\label{K3:l3}
For any $\u\in\fG_n$ and $B\in\B$,
$$
   \underset {d\sim\msD}{\Prob}\{\u(d)\in\U_B\} 
  = \underset{d\sim\msD}{\Prob}\{d\in\U_B\}
  =\frac{1}{2^n}.
$$
\end{lemma}

\proof
The equality between probabilities follows from~(i) 
in Definition~\ref{K3:def:invariance}. 
Therefore, by Corollary~\ref{cor:K3_partition} 
and Definition~\ref{K3:def:invariance}(ii),
the probability of each of 
them is $2^{-n}$. 
\eproof

The following lemma tells us
that, for all $B\in\B$, the random variable 
$h_B(d)$ is independent of the event ``$d\in\U_B$."

\begin{lemma}\label{K3:l4}
For all measurable $g:\R\rightarrow\R$ and $B\in\B$, 
$$
   \underset {d\sim\msD}{\E}
    \big(g\big(h_B(d)\big)\big| \,d\in\U_B\big) =
   \underset {d\sim\msD}{\E}
    \big(g\big(h_B(d)\big)\big).
$$
\end{lemma}

\proof
From the definition of conditional expectation 
and Lemma~\ref{K3:l3} we have 
\begin{equation}\label{K3:eq:cond_ex}
   \underset {d\sim\msD}{\E}
    \big(g(h_B(d))\big| \,d\in\U_B\big) =
   \frac{\displaystyle\int_{d\in\U_B}g(h_B(d))f(d)}
  {\underset{d\sim\msD}{\Prob}\{d\in\U_B\}}
=\;  2^n \int_{d\in\D}\uno_B(d) g(h_B(d))f(d)
\end{equation}
where $\uno_B$ denotes the indicator function of $\U_B$. 
Now, for any $\u\in\fG_n$, the map $d\mapsto \u(d)$ 
is a linear isometry on $\D$. Therefore
$$
  \int_{d\in\D}\uno_B(d) g(h_B(d))f(d)=
  \int_{d\in\D}\uno_B(\u(d)) g(h_B(\u(d)))f(\u(d)).
$$
Using that 
$h_B(d)=h_B(\u(d))$ (by Lemma~\ref{K3:l1}) 
and $f(d) = f(\u(d))$ (by the $\fG_n$-invariance of $\msD$), 
it follows that 
\begin{align*}
   \underset {d\sim\msD}{\E}
    \big(g(h_B(d))&\,\big| \,d\in\U_B\big) 
= 2^n \int_{d\in\D}\uno_B(d)g(h_B(d))f(d)\\
= &\sum_{\u\in\fG_n} \int_{d\in\D}\uno_B(\u(d)) 
     g(h_B(\u(d))) f(\u(d))\\
= &\sum_{\u\in\fG_n} \int_{d\in\D}\uno_B(\u(d)) 
     g(h_B(d)) f(d)\\
= &\int_{d\in\D} g(h_B(d))f(d)\;=\; 
    \underset {d\sim\msD}{\E}
    \big(g\big(h_B(d)\big)\big),
\end{align*}
the last line by Corollary~\ref{cor:K3_partition}.
\eproof

Let $B^*=\{1,2,\ldots,m\}$.

\begin{lemma}\label{K3:l5}
For all measurable $g:\R\rightarrow\R$
$$
  \underset {d\sim N(0,\Id)}{\E}(g(\varrho(d))\mid d\in\U)
  =\underset{d\sim N(0,\Id)}{\E}\big(g\big(h_{B^*}(d)\big)\big).
$$
\end{lemma}

\proof
Let $\varphi$ be the probability density function of $N(0,\Id)$.
\begin{equation}\label{K3:eq12}
  \underset{d\sim N(0,\Id)}{\E}(g(\varrho(d))|d\in\U)
  =\frac{\int_{d\in\U}g(\varrho(d))\varphi(d)\sfd(d)}
  {\underset{d\sim N(0,\Id)}{\Prob}\{d\in\U\}}.
\end{equation}
Since $d$ is Gaussian, the probability
that $d$ has two optimal bases is $0$. Using this and 
Lemma \ref{K3:l3} we see that 
\begin{equation}\label{K3:eq13}
  \underset{d\sim N(0,\Id)}{\Prob}\{d\in\U\}
 =\sum_{B\in\B}\underset{d\sim N(0,\Id)}{\Prob}\{d\in\U_B\}
 = \sum_{B\in\B}\frac{1}{2^n}
       =\left(\begin{array}{c}n\\m\end{array}\right)
       \left(\frac{1}{2^n}\right).
\end{equation}
Combining Equations (\ref{K3:eq12}) and (\ref{K3:eq13}), we have
\begin{eqnarray*}
   \left(\begin{array}{c}n\\m\end{array}\right)
  \left(\frac{1}{2^n}\right)\underset{d\sim N(0,\Id)}{\E}(g(\varrho(d))|d\in\U)
 &=&\int_{d\in\U}g(\varrho(d))\varphi(d)\sfd(d)\\
 &=&\sum_{B\in\B}\int_{d\in\U_B}g(\varrho(d))\varphi(d)\sfd(d)
\end{eqnarray*}
the last since the probability that $d$ has two optimal bases is
$0$. Using now that the entries of $d$ are i.i.d. 
and Theorem~\ref{K3:the1} we obtain
\begin{eqnarray*}
  \left(\frac{1}{2^n}\right)\underset{d\sim N(0,\Id)}{\E}
  (g(\varrho(d))\mid d\in\U)
  &=&\int_{d\in\U_{B^*}}g(\varrho(d))\varphi(d)\sfd(d)\\
  &=&\int_{d\in\U_{B^*}}g\big(h_{B^*}(d)\big)\varphi(d)\sfd(d).
\end{eqnarray*}
Therefore, by Lemma \ref{K3:l3} with $B=B^*$,
\begin{equation*}
  \underset{d\sim N(0,\Id)}{\Prob}\{d\in\U_{B^*}\}\,
 \underset{d\sim N(0,\Id)}{\E}(g(\varrho(d))\mid d\in\U)
 =\int_{d\in\U_{B^*}}
  g\big(h_{B^*}(d)\big)\varphi(d)\sfd(d).
\end{equation*}
We conclude since, by the definition of conditional expectation 
and Lemma \ref{K3:l4},
\begin{align}
   \underset {d\sim N(0,\Id)}{\E}(g(\varrho(d))\mid d\in\U) 
  \;&=\;\underset{d\sim N(0,\Id)}{\E} \big(g\big(h_{B^*}(d)
   \big)\mid d\in\U_{B^*}\big)\notag\\
  &=\;\underset{d\sim N(0,\Id)}{\E}
    \big(g\big(h_{B^*}(d)\big)\big).\tag*{\qed}
\end{align}

The following is Lemma~11 in~\cite{ChC02}. 

\begin{lemma}\label{K3:l6}
For the $\|\ \|_{12}$ in the definition of 
$\rho_\Sing$ we have 
$$
  \underset{S\sim N(0,\Id)}{\E} \left(\sqrt{\frac{1}{\rho_\Sing(S)}}\right) 
   \leq 2m^{5/4}
$$
where $N(0,\Id)$ is the Gaussian distribution in the set of 
$m\times m$ real matrices. 
\end{lemma}

\begin{lemma}\label{K3:l7}
Let $B\in\B$ fixed. Then, for the $\|\ \|_{12}$ 
in the definition of $\rho_\Sing$ we have 
$$
   \underset{d\sim N(0,\Id)}{\E}\left(\sqrt\frac{1}
  {h_{B}(d)}\right)\leq2(m+1)^{5/4}(n+1).
$$
\end{lemma}

\proof
For any fixed $d\in\D$,
$$
\sum_{S\in\S_{B}}\sqrt{\frac{1}{\rho_\Sing(S)}}
\,\geq\,\max_{S\in\S_{B}}\sqrt\frac{1}{\rho_\Sing(S)}
\,=\,\sqrt\frac{1}{h_{B}(d)}.
$$
Take average on both sides,
\begin{align}
   \underset{d\sim\msD}{\E}\left(\sqrt\frac{1}{h_{B}(d)}\right)
  \leq&\underset{d\sim\msD}{\E}\left(\sum_{S\in\S_{B}}
  \sqrt{\frac{1}{\rho_\Sing(S)}}\right)
  \,\leq\,\sum_{S\in\S_{B}}\underset{d\sim\msD}{\E}
   \left(\sqrt{\frac{1}{\rho_\Sing(S)}}\right)\notag\\
 \leq&\sum_{S\in\S_{B}}2(m+1)^{5/4}
  \qquad\mbox{by Lemma \ref{K3:l6}}\notag\\
\leq&\;2(m+1)^{5/4}(n+1).\tag*{\qed}
\end{align}
\medskip

The following lemma is proved as Lemma \ref{K3:l4}. 

\begin{lemma}\label{K3:l8}
For all $r,s\geq 1$ we have
\begin{equation}\tag*{\qed}
  \underset {d\sim\msD}{\E}\left(\|d\|_{rs}\mid d\in\U\right) 
  =\underset{d\sim\msD}{\E}\left(\|d\|_{rs}\right).
\end{equation}
\end{lemma}

\begin{lemma}\label{K3:l9}
We have
$$
  \underset {d\sim N(0,\Id)}{\E}\left(\|d\|_{12}\right)\leq 6\sqrt{n+1}.
$$
\end{lemma}

\proof
Recall that $\|d\|_{12}=\|M_d\|_{12}$. It is well known that 
$\|M_d\|_{12}\leq \|M_d\|$ where the latter is spectral 
norm. The statement now follows from the fact that, 
for a random Gaussian $A\in\R^{(m+1)\times (n+1)}$ 
we have $\E(\|A\|)\leq 6\sqrt{n+1}$~\cite[Lemma~2.4]{BuCu:10}.
\eproof 

\proofof{Theorem~\ref{K3:th:mainK}}
By Jensen's inequality and Lemma~\ref{K3:l9},
\begin{equation}\label{K3:eq14}
\underset{d}{\E}\left(\ln\|d\|_{12}\right) 
\leq\ln \underset{d}{\E}\left(\|d\|_{12}\right)
\leq \frac12\ln (n+1) +\ln 6.
\end{equation}
In addition, using now Lemma~\ref{K3:l7},
\begin{eqnarray}\label{K3:eq15}
   \underset{d}{\E}\left(\ln\left(h_{B^*}(d)\right)\right)
&=&
   -2\underset{d}{\E}\left(\ln\sqrt\frac{1}{h_{B^*}(d)}\right)
\,\geq\,
   -2\ln\underset{d}{\E}\left(\sqrt\frac{1}{h_{B^*}(d)}\right)\notag\\
&\geq& -\ln(2(m+1)^\frac{5}{4}(n+1)). 
\end{eqnarray}
By the definition of $\K(d)$ and Lemmas \ref{K3:l8} and \ref{K3:l5},
\begin{eqnarray}\label{K3:eq16}
\underset {d}{\E}\left(\left.\ln\K(d)\right|\,d\in\U\right)
&=&\underset{d}{\E}\left(\left.\ln\|d\|_{12}\right|\,d\in\U\right)
-\underset{d}{\E}\left(\left.\ln\varrho(d)\right|\,d\in\U\right)\notag\\
&=&\underset{d}{\E}\left(\ln\|d\|_{12}\right)
    -\underset{d}{\E}\left(\ln\left(h_{B^*}(d)\right)\right).
\end{eqnarray}
Combining Equations (\ref{K3:eq14}), (\ref{K3:eq15}), and 
(\ref{K3:eq16}), 
the proof is done.
\eproof

{\small

}

\end{document}